\documentclass[final]{amsart}
\usepackage{showkeys}
\usepackage{amsmath,amsfonts,amssymb,latexsym,epic,eepic}
\usepackage[dvips]{graphics,epsfig,color}
\usepackage{pstricks}

\newtheorem{theorem}{Theorem}[section]
\newtheorem{lemma}[theorem]{Lemma}
\theoremstyle{definition}

\theoremstyle{remark}
\newtheorem{remark}[theorem]{Remark}
\numberwithin{equation}{section}

\newcommand{\crou}{\phi} 
\newcommand{\ie}{\emph{i.e.\/}}
\newcommand{\eg}{\emph{e.g.\/}}
\newcommand{\dx}{\, {\rm d}\bfx}
\newcommand{\dedge}{\, {\rm d} \gamma}
\newcommand{\ds}{\, {\rm d}s}

\newcommand{\mach}{{\rm Ma}}
\newcommand{\machd}{A}
\newcommand{\xN}{\mathbb{N}}
\newcommand{\xR}{\mathbb{R}}
 
\newcommand{\xLtwo}{{\rm L}^2} 
\newcommand{\xLinfty}{{\rm L}^\infty} 
\newcommand{\xHone}{{\rm H}^{1}}
\newcommand{\xHn}[1]{{\rm H}^#1}
\newcommand{\xCinfty}{{\rm C}^{\infty}} 
\newcommand{\edge}{\sigma}

\newcommand{\edges}{{\mathcal E}}
\newcommand{\edgesint}{{\mathcal E}_{{\rm int}}}
\newcommand{\edgesext}{{\mathcal E}_{{\rm ext}}}
\newcommand{\mesh}{{\mathcal T}}
\newcommand{\vs}{\mathrm{v}_{\edge,K}}
\newcommand{\vsp}{\mathrm{v}_{\edge,K}^+}
\newcommand{\vsm}{\mathrm{v}_{\edge,K}^-}

\newcommand{\normLd}[1]{\hspace{.2em}|\hspace{-.1em}| #1 |\hspace{-.1em}|_{\xLtwo(\Omega)}\hspace{.2em}}
\newcommand{\normLdom}[2]{\hspace{.2em}|\hspace{-.1em}| #1 |\hspace{-.1em}|_{\xLtwo(#2)}\hspace{.2em}}
\newcommand{\normLdomd}[2]{\hspace{.2em}|\hspace{-.1em}| #1 |\hspace{-.1em}|_{\xLtwo(#2)}^2\hspace{.2em}}

\newcommand{\normLddR}[1]{\hspace{.2em}|\hspace{-.1em}| #1 |\hspace{-.1em}|_{\xLtwo(\xR^d)}^2\hspace{.2em}}

\newcommand{\snormH}[1]{\hspace{.2em}| #1 |_{\xHone(\Omega)}\hspace{.2em}}
\newcommand{\normHb}[1]{\hspace{.2em}|\hspace{-.1em}| #1 |\hspace{-.1em}|_{1,b}\hspace{.2em}}
\newcommand{\normHbd}[1]{\hspace{.2em}|\hspace{-.1em}| #1 |\hspace{-.1em}|_{1,b}^2\hspace{.2em}}
\newcommand{\snormdisc}[1]{\hspace{.2em}| #1 |_{{\mesh,\beta}}\hspace{.2em}}
\newcommand{\snormdiscd}[1]{\hspace{.2em}| #1 |_{{\mesh\hspace{-.2em},\beta}}^2\hspace{.2em}}
\newcommand{\snormHdeuxdom}[2]{\hspace{.2em}| #1 |_{\xHn{2}(#2)}\hspace{.2em}}

\newcommand{\bff}{{\boldsymbol f}}
\newcommand{\bfn}{{\boldsymbol n}}
\newcommand{\bfq}{{\boldsymbol q}}
\newcommand{\bft}{{\boldsymbol t}}
\newcommand{\bfu}{{\boldsymbol u}}
\newcommand{\bfv}{{\boldsymbol v}}
\newcommand{\bfWh}{{\boldsymbol W}_{\hspace{-.2em}h}}
\newcommand{\bfx}{{\boldsymbol x}}
\newcommand{\bfy}{{\boldsymbol y}}
\newcommand{\grad}{{\boldsymbol \nabla}}
\newcommand{\bfeta}{{\boldsymbol \eta}}
\newcommand{\bfpsi}{{\boldsymbol \psi}}
\newcommand{\dive}{{\rm div}}

\begin{document}

\title[A FE-FV scheme for the isothermal compressible {S}tokes problem]
{A convergent Finite Element-Finite Volume scheme for the compressible {S}tokes problem \\ Part I -- the isothermal case}

\author{T. Gallou\"et}
\address{Universit\'e de Provence, France}
\email{gallouet@cmi.univ-mrs.fr}

\author{R. Herbin}
\address{Universit\'e de Provence, France}
\email{herbin@cmi.univ-mrs.fr}

\author{J.-C. Latch\'e}
\address{Institut de Radioprotection et de S\^{u}ret\'{e} Nucl\'{e}aire (IRSN)}
\email{jean-claude.latche@irsn.fr}

%
%
%
%
%
%


\begin{abstract}
In this paper, we propose a discretization for the (nonlinearized) compressible Stokes problem with a linear equation of state $\rho=p$, based on Crouzeix-Raviart elements.
The approximation of the momentum balance is obtained by usual finite element techniques.
Since the pressure is piecewise constant, the discrete mass balance takes the form of a finite volume scheme, in which we introduce an upwinding of the density, together with two additional stabilization terms.
We prove {\em a priori} estimates for the discrete solution, which yields its existence by a topological degree argument, and then the convergence of the scheme to a solution of the continuous problem.
\end{abstract}

\maketitle


\section{Introduction}

The problem addressed in this paper is the system of the so-called barotropic compressible Stokes equations, which reads:
\begin{subequations} \begin{align}
& \displaystyle \label{pbcont1}
- \Delta \bfu + \grad p = \bff
\\[1ex] & \displaystyle \label{pbcont2}
\dive (\rho\, \bfu)=0 
\\[1ex] & \displaystyle \label{pbcont3}
\rho=\varrho(p)
\end{align} \label{pbcont} \end{subequations}
where $\rho$, $\bfu$ and $p$ stand for the density, velocity and pressure in the flow, respectively, and $\bff$ is a forcing term.
The function $\varrho$ is the equation of state used for the modelling of the particular flow at hand, which may be the actual equation of state of the fluid or may result from assumptions concerning the flow.
Here, we only consider the following relation, which corresponds to an isothermal flow of a perfect gas:
\begin{equation}
\varrho(p) = \machd \ p 
\label{eos}\end{equation}
where $\machd$ is a positive constant.
Since the sound velocity is defined by $c^2=dp/d\rho$, $\machd=\mach^2/V^2$, where $\mach$ is the Mach number and $V$ is a characteristic velocity of the flow.
This system of equations is posed over $\Omega$, a bounded domain of $\xR^d$, $d\leq 3$ supposed to be polygonal ($d=2$) or polyhedral ($d=3$).
It is supplemented by homogeneous boundary conditions for $\bfu$, and by prescribing the total mass $M$ of the fluid:
\begin{equation}
\int_\Omega \rho \dx=M
\label{totmass}\end{equation}
where $M$ is a positive real number.

\medskip
In this paper, we study a numerical scheme for the solution of this problem, which combines low order finite element and finite volume techniques, and is very close to a scheme which was proposed for the solution of barotropic Navier-Stokes equations in \cite{gal-08-unc} and further extended to two-phase flows in \cite{gas-08-ent}; the resulting code is today currently used at the French {\it Institut de Radioprotection et de S\^uret\'e Nucl\'eaire} (IRSN) for "real-life" studies in the nuclear safety field.
Up to now, stability (in the sense of conservation of the entropy) is known for these schemes, and numerical experiments show convergence rates close to one in natural energy norms.
Our goal is now to prove their convergence.
This work is the first one in this direction, and we address here the probably simplest toy problem, restricting ourselves to the steady case, to creeping flows (\ie\ omitting the convection term in the momentum balance equation) and to a linear equation of state.
The extension to laws where $\rho$ varies linearly with $p^{1/\gamma}$, where $\gamma>1$ is a coefficient which is specific to the considered fluid, which are typically obtained for isentropic flows of perfect gases, is the object of a further paper (part II of the present one); the additional difficulty posed by this further study is to prove the strong convergence for the density, which necessitates to adapt P.L. Lions' ``effective pressure trick" \cite{pll-98-mat}  at the discrete level.
Finally, for the sake of simplicity, we use here a simplified form of the diffusion term ($- \Delta \bfu$) but it is clear from the subsequent developements that the presented theory holds for any linear elliptic operator (and in particular for the usual form of the viscous term for compressible constant viscosities flows).

\medskip
The finite element - finite volume discretization which is chosen here is motivated by the fact that we wish the approximate density to be positive, as in the continuous model,  in order to be compatible with the physics.
Moreover, the proof of convergence of a numerical approximation to \eqref{pbcont} requires estimates on both velocity and pressure or density, and the density positivity is very useful to obtain these estimates.
A classical way to ensure positivity is  to use a finite volume upwinding technique in the discretization of the term $\dive (\rho u)$.
This technique is easily set up if the discrete velocities are located on the edges and densities and pressures at the cell centres, which is the reason why we choose the Crouzeix-Raviart finite elements for the velocities and cell centred finite volumes for the densities.

\medskip
This paper is organized as follows.
The discretization is first described (section \ref{sec:disc}), and we prove an $\xLtwo$ compactness result for sequences of Crouzeix-Raviart functions with bounded broken $\xHone$ semi-norm (section \ref{sec:comp}).
Then the proposed scheme is given (section \ref{sec:scheme}), and the above-mentioned compactness result yields the convergence of (sub-)sequences of discrete solutions to a limit, thanks to {\em a priori} estimates which are given in section \ref{sec:apriori}.
Finally, this limit is shown to be a solution to the continuous problem in section \ref{sec:conv}.

\medskip
To our knowledge, this convergence proof is the first one for the genuine (nonlinear) compressible Stokes problem; a linearized version of this system is adressed in previous works \cite{kel-96-fin,kel-97-pen,kwe-00-opt,kwe-03-opt,boc-06-ana}.


\section{Discrete spaces and relevant lemmata} \label{sec:disc}

Let $\mesh$ be a decomposition of the domain $\Omega$ in simplices. 
By $\edges(K)$, we denote the set of the edges ($d=2$) or faces ($d=3$) $\edge$ of the element $K \in \mesh$; for short, each edge or face will be called an edge hereafter.
The set of all edges of the mesh is denoted by $\edges$; the set of edges included in the boundary of $\Omega$ is denoted by $\edgesext$ and the set of internal ones (\ie\ $\edges \setminus \edgesext$) is denoted by $\edgesint$.
The decomposition $\mesh$ is supposed to be regular in the usual sense of the finite element literature (\eg\ \cite{cia-91-bas}), and, in particular, $\mesh$ satisfies the following properties: $ \bar\Omega=\bigcup_{K\in \mesh} \bar K$; if $K,\,L \in \mesh,$ then $\bar K \cap \bar L=\emptyset$, $\bar K \cap \bar L=\emptyset$ is a vertex of the mesh or $\bar K\cap \bar L$ is a common edge of $K$ and $L$, which is denoted by $K|L$.
For each internal edge of the mesh $\edge=K|L$, $\bfn_{KL}$ stands for the normal vector of $\edge$, oriented from $K$ to $L$ (so that $\bfn_{KL}=-\bfn_{LK}$).
By $|K|$ and $|\edge|$ we denote the measure, respectively, of the element $K$ and of the edge $\edge$, and $h_K$ and $h_\edge$ stand for the diameter of $K$ and $\edge$, respectively.
We measure the regularity of the mesh through the parameter $\theta$ defined by:
\begin{equation}
\theta=\inf \ \lbrace \frac{\xi_K}{h_K};\ K \in \mesh \rbrace
\cup \lbrace  \frac{h_L}{h_K}, \frac{h_K}{h_L};\ \edge=K|L \in \edgesint \rbrace,
\label{regul}\end{equation}
where $\xi_K$ stands for the diameter of the largest ball included in $K$.
Note that the following inequality holds:
\begin{equation}
h_\edge\, |\edge| \leq 2\ \theta^{-d} \, |K|, \, \forall K \in\mesh, \, \forall \edge \in \edges(K).
\label{ineqh}
\end{equation}
Indeed, this relation is derived by noting that $h_\edge\, |\edge| \leq h_K^d$ and $|K| \geq c\,\xi_K^d$ with $c=\pi/4$ in 2D and $c=\pi/6$ in 3D;
it will be used throughout this paper.
Finally, as usual, we denote by $h$ the quantity $\max_{K\in\mesh} h_K$.

\smallskip
The space discretization relies on the Crouzeix-Raviart element (see \cite{cro-73-con} for the seminal paper and, for instance, \cite[p. 199--201]{ern-04-the} for a synthetic presentation).
The reference element is the unit $d$-simplex and the discrete functional space is the space $P_1$ of affine polynomials.
The degrees of freedom are determined by the following set of nodal functionals:
\begin{equation}
\displaystyle \left\{F_\edge,\ \edge \in \edges(K) \right\}, 
\qquad F_\edge(v)=|\edge|^{-1}\int_{\edge} v \dedge.
\label{vdof}\end{equation}
The mapping from the reference element to the actual one is the standard affine mapping.
Finally, the continuity of the average value of the discrete functions (\ie, for any function $v$, $F_\sigma(v)$) across each face of the mesh is required, thus the discrete space $V_{h}$ is defined as follows:
\begin{equation}
\begin{array}{ll} \displaystyle
V_h =
& \displaystyle
\lbrace
\, v\in L^{2}(\Omega)\,:\, v|_K \in P_1(K),\,\forall K\in \mesh;\ F_\edge (v) \mbox{ continuous} \\[1ex] & \displaystyle\mbox{  across}

\mbox{ each edge } \sigma \in {\edgesint};\ \ F_\edge(v)=0,\ \forall \edge \in \edgesext \rbrace.
\end{array}
\label{crdef}\end{equation}
The space of approximation for the velocity is the space $\bfWh$ of vector valued functions each component of which belongs to $V_h$: $\bfWh=(V_h)^d$.
The pressure is approximated by the space $L_{h}$ of piecewise constant functions:
\[
L_h=\left\{q\in L^{2}(\Omega)\,:\, q|_K=\mbox{ constant},\,\forall K\in \mesh\right\}.
\]
Since only the continuity of the integral over each edge of the mesh is imposed, the functions of $V_h$ are discontinuous through each edge; the discretization is thus nonconforming in $H^1(\Omega)^d$.
We then define, for $1 \leq i \leq d$ and $v\in V_h$, $\partial_{h,i}\, v$ as the function of $\xLtwo(\Omega)$ which is equal to the (piecewise constant) derivative of $v$ with respect to the $i^{th}$ space variable almost everywhere.
This notation allows to define the discrete gradient, denoted by $\grad_h$, for both scalar and vector valued discrete functions and the discrete divergence of vector valued discrete functions, denoted by $\dive_h$.

\medskip
The Crouzeix-Raviart pair of approximation spaces for the velocity and the pressure is \textit{inf-sup} stable, in the usual sense for "piecewise $\xHone$" discrete velocities, \ie\ there exists $c_{\rm i}>0$ independent of the mesh such that:
\[
\forall q \in L_h, \qquad \sup_{\bfv \in \bfWh} \frac{\displaystyle \sum_{K\in \mesh} \int_K q \ \dive \bfv \dx}{\normHb{\bfv}} =
\sup_{\bfv \in \bfWh} \frac{\displaystyle \int_\Omega q \ \dive_h \bfv \dx}{\normHb{\bfv}}
\geq c_{\rm i} \normLd{q- q_{\rm m}},
\]
where $q_{\rm m}$ is the mean value of $q$ over $\Omega$ and $\normHb{\cdot}$ stands for the broken Sobolev $\xHone$ semi-norm, which is defined for any function $v \in V_h$ or  $v \in \bfWh$ by:
\[
\normHbd{v}=\sum_{K\in \mesh} \int_K |\grad v |^2 \dx=\int_\Omega | \grad_h v |^2 \dx.
\]
This broken Sobolev semi-norm is known to control the $\xLtwo$ norm by an extended Poincar\'e inequality \cite[proposition 4.13]{tem-77-nav}, in the sense that for any function $v \in V_h$, $\normHb{v} \leq c \normLd{v}$ where the real number $c$ only depends on the computational domain.

We also define a discrete semi-norm on $L_h$, similar to the usual finite volume discrete $\xHone$ semi-norm, weighted by a mesh-dependent coefficient:
\[
\forall q \in L_h, \qquad
\snormdiscd{q}=
\sum_{\stackrel{\scriptstyle \edge \in \edgesint,}{\scriptstyle \edge=K|L}} (h_K+h_L)^\beta\ \frac{|\edge|}{h_\edge}\ (q_K-q_L)^2.
\]

\medskip
From the definition \eqref{vdof}, each velocity degree of freedom can be associated to an element edge.
Consequently, the velocity degrees of freedom are indexed by the number of the component and the associated edge, thus the set of velocity degrees of freedom reads:
\[
\lbrace u_{\edge, i},\ \edge \in \edgesint,\ 1 \leq i \leq d \rbrace.
\]

We denote by $\crou_\edge$  the usual Crouzeix-Raviart shape function associated to $\edge$, \ie\ the scalar function of $V_h$ such that $F_\edge(\crou_\edge)=1$ and $F_{\edge'}(\crou_\edge)=0,\ \forall \edge' \in \edges \setminus \{ \edge \}$.

\medskip
Similarly, each degree of freedom for the pressure is associated to a cell $K$, and the set of pressure degrees of freedom is denoted by
$\lbrace p_K,\ K \in \mesh \rbrace$.

\bigskip
Finally, we define by $r_h$ the following interpolation operator:
\begin{equation}
\begin{array}{l|lcl}
r_h : \qquad
&
\xHone_0(\Omega)  & \longrightarrow & V_h
\\ &
v & \mapsto & \displaystyle 
r_h v = \sum_{\edge \in \edges} F_\edge(v)\,\crou_\edge =\sum_{\edge \in \edges} |\edge|^{-1} \left(\int_{\edge} v \dedge \right) \,\crou_\edge.
\end{array}
\label{def_rh}\end{equation}
This operator naturally extends to vector-valued functions (\ie\ to perform the interpolation from $\xHone_0(\Omega)^d$ to $\bfWh$), and we keep the same notation $r_h$ for both the scalar and vector case.
The properties of $r_h$ are gathered in the following lemma.
They are proven in \cite{cro-73-con}.

\begin{lemma}
Let $\theta_0 >0$ and let $\mesh$ be a triangulation of the computational domain $\Omega$ such that $\theta \geq \theta_0$, where $\theta$ is defined by \eqref{regul}.
The interpolation operator $r_h$ enjoys the following properties:
\begin{enumerate}
\item Preservation of the divergence:
\begin{equation}
\forall \bfv \in \xHone_0(\Omega)^d,\ \forall q \in L_h, \qquad \int_\Omega q\ \dive_h (r_h \bfv) \dx = \int_\Omega q\ \dive \bfv \dx.
\label{interp_1}\end{equation}
\item Stability:
\begin{equation}
\forall v \in \xHone_0(\Omega),\qquad \normHb{r_h v} \leq c_1(\theta_0) \snormH{v}.
\label{interp_2}\end{equation}
\item Approximation properties:
\begin{equation}
\begin{array}{l}\displaystyle
\forall v \in \xHn{}2(\Omega) \cap \xHone_0(\Omega),\ \forall K \in \mesh, 
\\[1ex] \displaystyle \hspace{5ex}
\normLdom{v-r_h v}{K} + h_K \normLdom{\grad_h(v-r_h v)}{K} \leq c_2(\theta_0)\, h_K^2 \snormHdeuxdom{v}{K}.
\end{array}
\label{interp_3}\end{equation}
\end{enumerate}
\label{interpolation}\end{lemma}

In both above inequalities, the notation $c_i(\theta_0)$ means that the real number $c_i$ only depends on $\theta_0$, and, in particular, not on the parameter $h$ characterizing the size of the cells; this notation will be kept throughout the paper.

\medskip
The following lemma is known (\eg \ \cite[lemma 3.32]{ern-04-the}); we give its (elementary) proof for the sake of completeness.
\begin{lemma}
Let $\theta_0 >0$ and let $\mesh$ be a triangulation of the computational domain $\Omega$ such that $\theta \geq \theta_0$, where $\theta$ is defined by \eqref{regul}, and $V_h$ be the space of Crouzeix-Raviart discrete functions associated to $\mesh$, as defined by \eqref{crdef}.
Then there exists a real number $c(\theta_0)$ such that the following bound holds for any $v \in V_h$:
\[
\sum_{\edge \in \edges} \frac 1 {h_\edge} \int_\edge [v]^2 \dedge \leq c(\theta_0) \normHbd{v},
\]
where, on any $\edge \in \edgesint$, $[v]$ stands for the jump of $v$ across $\edge$ and, on any $\edge \in \edgesext$, $[v]=v$.
\label{est_jumps}\end{lemma}

\begin{proof}
For any control volume $K$ of the mesh, we denote by $(\grad v)_K$ the (constant) gradient of the restriction of $v$ to $K$.
With this notation, using the continuity of $v$ across $\edge$ at the mass center $\bfx_\edge$ of any internal edge and the fact that $v$ vanishes at the mass center $\bfx_\edge$ of any external edge, we get:
\[
\begin{array}{l} \displaystyle
\sum_{\edge \in \edges} \frac 1 {h_\edge} \int_\edge [v]^2 \dedge=
\sum_{\stackrel{\scriptstyle \edge \in \edgesint,}{\scriptstyle \edge=K|L}} 
   \frac 1 {h_\edge} \int_\edge \left( ((\grad v)_K-(\grad v)_L) \cdot (\bfx-\bfx_\edge)\right)^2 \dedge
\\ \displaystyle \hspace{30ex}
+ \sum_{\stackrel{\scriptstyle \edge \in \edgesext,}{\scriptstyle \edge \in \edges(K)}}
\frac 1 {h_\edge} \int_\edge \left( (\grad v)_K \cdot (\bfx-\bfx_\edge)\right)^2 \dedge.
\end{array}
\]
We thus have:
\[
\sum_{\edge \in \edges} \frac 1 {h_\edge} \int_\edge [v]^2 \dedge \leq
2 \sum_{\stackrel{\scriptstyle \edge \in \edgesint,}{\scriptstyle \edge=K|L}} h_\edge \, |\edge|\ (|(\grad v)_K|^2 + |(\grad v)_L|^2)
+ \sum_{\stackrel{\scriptstyle \edge \in \edgesext,}{\scriptstyle \edge \in \edges(K)}} h_\edge \, |\edge|\ |(\grad v)_K|^2.
\]
and the result follows by regularity of the mesh.
\end{proof}

The proof of the following trace lemma can be found in \cite[section 3]{ver-99-err}.

\begin{lemma}\label{trace}
Let $\mesh$ be a given triangulation of $\Omega$ and $K$ be a control volume of $\mesh$, $h_K$ its diameter and $\edge$ one of its edges.
Let $v$ be a function of $\xHone(K)$.
Then the following inequality holds:
\[
\normLdom{v}{\edge}\leq \left( d\ \frac{|\edge|}{|K|} \right)^{1/2}
\ \left ( \normLdom{v}{K} + h_K \normLdom{\nabla v}{K} \right ).
\]
\end{lemma}
We will also need the following Poincar\'e ineqality:
\begin{equation}
\forall K \in \mesh,\ \forall v \in \xHone(K), \qquad \normLdom{v - v_{{\rm m},K}}{K} \leq \frac 1 \pi \ h_K \normLdom{\nabla v}{K}.
\label{poincare}\end{equation}
where $v_{{\rm m},K}$ stands for the mean value of $v$ over $K$.
This relation is proven for any convex domain in \cite{pay-60-opt}.

\bigskip
We are now in position to prove the following technical result.

\begin{lemma}
Let $\theta_0 >0$ and let $\mesh$ be a triangulation of the computational domain $\Omega$ such that $\theta \geq \theta_0$, where $\theta$ is defined by \eqref{regul}; let $(a_\edge)_{\edge\in \edgesint}$ be a family of real numbers such that $\forall \edge \in \edgesint,\ a_\edge\leq 1$ and let $v$ be a function of the Crouzeix-Raviart space $V_h$ associated to $\mesh$.
Then the following bound holds:
\[
\sum_{\edge \in \edgesint} \left| \int_\edge a_\edge\ [v]\, f \dedge \right| \leq c(\theta_0) \, h \normHb{v} \snormH{f}, \, \forall f\in \xHone_0(\Omega).
\]
where the real number $c(\theta_0)$ only depends on $\theta_0$ and on the domain $\Omega$.
\label{jump_reg}\end{lemma}
\begin{proof}
Since the integral of the jump across any edge of the mesh of a function of $V_h$ is zero, we have, for any $\edge \in \edgesint$:
\[
\int_\edge a_\edge\ [v]\, f \dedge =\int_\edge a_\edge\ [v]\, (f-f_\edge) \dedge,
\]
where $f_\edge$ is any real number.
Using the Cauchy-Schwarz inequality, first in $\xLtwo(\edge)$ then in $\xR^{card(\edges)}$ we thus get:
\[
\begin{array}{ll} \displaystyle
\sum_{\edge \in \edgesint} \left| \int_\edge a_\edge\ [v]\, f \dedge \right|
& \displaystyle \leq
\sum_{\edge \in \edgesint} \left[ \int_\edge [v]^2 \dedge \right]^{1/2}\  \left[ \int_\edge (f-f_\edge)^2 \dedge \right]^{1/2}
\\ & \displaystyle \leq
\left[ \sum_{\edge \in \edgesint} \frac 1 {h_\edge} \int_\edge [v]^2 \dedge  \right]^{1/2}
\underbrace{\left[ \sum_{\edge \in \edgesint} h_\edge \int_\edge (f-f_\edge)^2 \dedge  \right]^{1/2}}_{\displaystyle T_1}.
\end{array}
\]
By lemma \ref{est_jumps}, the first term of the latter product is bounded by $c(\theta_0) \normHb{v}$.
For the second one, choosing arbitrarily one adjacent simplex to each edge and applying the above trace lemma \ref{trace}, we get:
\[
T_1^2 \leq \sum_{\stackrel{\scriptstyle \edge \in \edgesint}{\scriptstyle (\edge \in \edges(K))}} 2 d\, h_\edge
\ \frac{|\edge|}{|K|}
\ \left (\normLdomd{f-f_\edge}{K} + h_K^2 \normLdomd{\nabla f}{K} \right).
\]
Choosing for $f_\edge$ the mean value of $f$ on $K$ and using \eqref{poincare}, we thus get:
\[
T_1^2 \leq \sum_{\stackrel{\scriptstyle \edge \in \edgesint}{\scriptstyle (\edge \in \edges(K))}} 
2d\, (1 + \frac 1 {\pi^2})\, h_\edge \, \frac{|\edge|}{|K|} h_K^2 \normLdomd{\nabla f}{K}.
\]
and the result follows by observing that the $\xHone$ semi-norm of $f$ on $K$ appears at most $(d+1)$ times in the summation and using the regularity of the mesh.
\end{proof}


\section{A compactness result}\label{sec:comp}
The aim of this section is to state and prove a compactness result for  $\normHb{\cdot}$  bounded sequences of discrete functions.
We begin by a preliminary lemma.
\begin{lemma}\label{ex_dedge}
Let $\theta_0 >0$ and let $\mesh$ be a triangulation of the computational domain $\Omega$ such that $\theta \geq \theta_0$, where $\theta$ is defined by \eqref{regul}; for $\edge \in \edges$, let $\chi_\edge$ be the function defined by:
\[
\begin{array}{l|lcl} \displaystyle
\chi_\edge\ : \hspace{1ex}
& \displaystyle
\xR^d \times \xR^d
& \longrightarrow &
\{0, \,1\}
\\[1ex] & \displaystyle
(\bfx,\bfy)
& \mapsto & \displaystyle
\chi_\edge(\bfx,\bfy)=1 \mbox{ if } [\bfx,\bfy] \cap \edge \neq \emptyset,\ \chi_\edge(\bfx,\bfy)=0 \mbox{ otherwise},
\end{array}
\]
where $\bfx$ and $\bfy$ are two points of $\xR^d$.
Then there exists a family of positive real numbers $(d_\edge)_{\edge \in \edges}$ such that:
\begin{enumerate}
\item for any $\edge \in \edges$, $d_\edge = c_1(\theta_0)\ h_\edge$,
\item for any points $\bfx$ and $\bfy$ of $\xR^d$ (possibly located outside $\Omega$), the following inequality holds:
\[
\sum_{\edge \in \edges} \chi_\edge(\bfx,\bfy)\  d_\edge  \leq c_2(\theta_0)\ (|\bfy-\bfx| + h)
\]
\end{enumerate}
\end{lemma}
\begin{proof}
We first deal with the two-dimensional case and with quasi-uniform meshes (\ie\ the bound we first prove blows up when $\max_{K\in\mesh} (h/h_K)$ tends to infinity).

\medskip
Let $\mesh$ be a triangulation of a two-dimensional domain $\Omega$, $K$ a triangular cell of $\mesh$ and $\edge$ an edge of $K$.
Without loss of generality, we suppose that $\edge$ is the segment $(0,h_\edge)\times {0}$ and we denote by $\xi_K$ the diameter of the largest ball included in $K$ and by $h_K$ the diameter of $K$.
We denote by $z_\edge$ the opposite vertex to $\edge$; the first coordinate of $z_\edge$ is necessarily lower than $h_K$ while its second coordinate is necessarily greater than $\xi_K$ (in the opposite case, no ball of diameter $\xi_K$ would be included in $K$).
It thus follows (see figure \ref{figmesh}):
\begin{enumerate}
\item that the rectangular domain $\omega_\edge=(h_\edge/3, 2 h_\edge/3) \times (0,h_\edge \xi_K /(12 h_K))$ is included in $K$,
\item that, if the similar construction is performed for another edge $\edge'$ of $K$ to obtain $\omega_{\edge'}$, $\omega_\edge$ and $\omega_\edge'$ do not intersect.
\end{enumerate}
We denote by $d_\edge$ the quantity $d_\edge=h_\edge \xi_K /(12 h_K)$.
We thus have $d_\edge \geq (\theta/12)\ h_\edge$, where $\theta$ is the parameter defined by \ref{regul}.

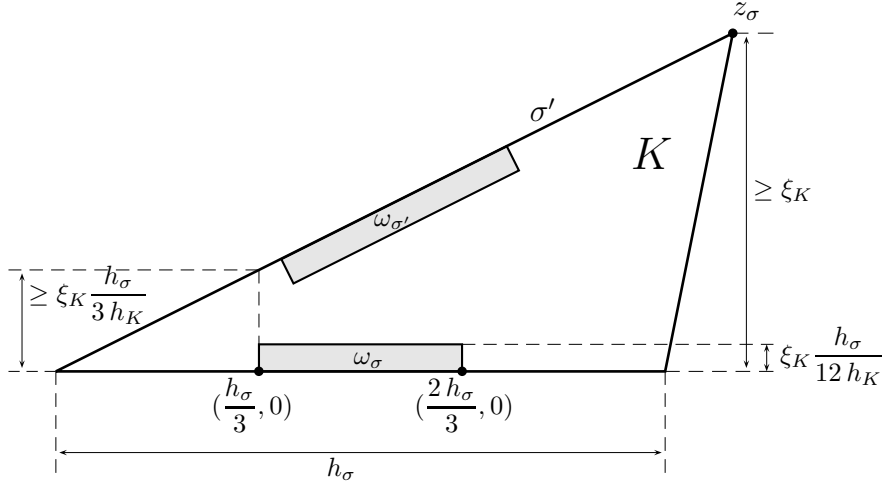
\begin{figure}[htb]
\begin{center}
\newgray{graydark}{.5}
\newgray{graymid}{.85}
\newgray{grayml}{.9}
\newrgbcolor{palombaggia}{0.6 0.65 1.}
\psset{unit=0.9cm}
\begin{pspicture}(0,0)(14,8)
\rput[bl](0,1){

   \rput[bl](9.5,4){\huge $K$}

   \pspolygon[fillstyle=solid, fillcolor=grayml](4.,1)(7,1)(7,1.4)(4,1.4)
   \rput[bl](5.4,1.05){$\omega_\edge$}
   \psline[linecolor=black, linewidth=0.3pt, linestyle=dashed]{-}(7,1.4)(11.6,1.4)
   \psline[linecolor=black, linewidth=0.3pt]{<->}(11.5,1)(11.5,1.4)
   \rput[bl](11.7,0.8){$\displaystyle \xi_K \frac {h_\edge}{12\,h_K}$}
   
   \pspolygon[fillstyle=solid, fillcolor=grayml](4.33,2.66)(7.66,4.33)(7.84,3.97)(4.51,2.3)
   \rput[bl](5.7,3.1){$\omega_{\edge'}$}
   \rput[bl](8,4.7){\Large $\edge'$}

   \psline[linecolor=black, linewidth=1pt]{-}(1,1)(10,1)(11,6)(1,1)
   \psdots*(11,6) \rput[bl](11,6.2){{\large $z_\edge$}}
   \psdots*(4,1) \rput[bl](3.3,0.1){$\displaystyle (\frac{h_\edge} 3,0)$}
   \psdots*(7,1) \rput[bl](6.3,0.1){$\displaystyle (\frac{2\,h_\edge} 3,0)$}
   
   \psline[linecolor=black, linewidth=0.3pt, linestyle=dashed]{-}(1,-0.5)(1,1)
   \psline[linecolor=black, linewidth=0.3pt, linestyle=dashed]{-}(10,-0.5)(10,1)
   \psline[linecolor=black, linewidth=0.3pt]{<->}(1.05,-0.2)(9.95,-0.2)
   \rput[bl](5,-0.6){$h_\edge$}
   
   \psline[linecolor=black, linewidth=0.3pt, linestyle=dashed]{-}(10,1)(11.6,1)
   \psline[linecolor=black, linewidth=0.3pt, linestyle=dashed]{-}(11,6)(11.6,6)
   \psline[linecolor=black, linewidth=0.3pt]{<->}(11.2,1.05)(11.2,5.95)
   \rput[bl](11.3,3.5){$\geq \xi_K$}
   
   \psline[linecolor=black, linewidth=0.3pt, linestyle=dashed]{-}(4,2.5)(0.3,2.5)é
   \psline[linecolor=black, linewidth=0.3pt, linestyle=dashed]{-}(1,1)(0.3,1)
   \psline[linecolor=black, linewidth=0.3pt, linestyle=dashed]{-}(4,1)(4,2.5)
   \psline[linecolor=black, linewidth=0.3pt]{<->}(0.5,1.05)(0.5,2.45)
   \rput[bl](0.6,1.7){$\displaystyle \geq \xi_K \frac {h_\edge}{3\,h_K}$}
}
\end{pspicture}
\caption{Notations for the control volume $K$}
\label{figmesh}
\end{center}
\end{figure}

\medskip
We now perform this construction for each edge $\edge$ of the mesh.
If $\edge\in \edgesext$, there is only one possible choice for $K$ (the adjacent cell to $\edge$); if $\edge \in \edgesint$, $\edge=K|L$, we choose either $K$ or $L$.
Let $\bfx$ and $\bfy$ be two points of $\xR^2$.
Let $\bft_{(\bfx,\bfy)}$ be the vector given by:
\[
\bft_{(\bfx,\bfy)}=\frac{\bfy-\bfx}{|\bfy-\bfx|}
\]
and $\bfn_{(\bfx,\bfy)}$ a normal vector to $\bft_{(\bfx,\bfy)}$.
We denote by $S_{(\bfx,\bfy)}$ the rectangle defined by:
\[
S_{(\bfx,\bfy)}=\left\{ \bfx + \xi_1\, \bft_{(\bfx,\bfy)} + \xi_2\, \bfn_{(\bfx,\bfy)},\ \xi_1 \in (-h,\ |\bfy-\bfx|+h)
,\ \xi_2 \in (-h,\ +h) \right\}
\]
For each edge intersected by the segment $[\bfx,\bfy]$ (\ie\ for each edge $\edge$ such that $\chi_\edge(\bfx,\bfy)=1$), the rectangle $\omega_\edge$ is included in $S_{(\bfx,\bfy)}$; thus, since these domains $\omega_\edge$ and $\omega_{\edge'}$ are disjoint:
\[
\sum_{\edge \in \edges} \chi_\edge(\bfx,\bfy)\  |\omega_\edge|  \leq |S_{(\bfx,\bfy)}|
\]
and thus:
\[
\sum_{\edge \in \edges} \chi_\edge(\bfx,\bfy)\  \frac 1 3 \ d_\edge\, h_\edge  \leq 2\ h \,(|\bfy-\bfx|+2h),
\]
which concludes the proof if $\max_{K\in\mesh} (h/h_K)$ is supposed to be bounded.

\medskip
The extension to the three-dimensional case is straightforward, since it only necessitates to adapt the construction of the domains $\omega_\edge$.
Finally, giving up the assumption that $\max_{K\in\mesh} (h/h_K)$ is bounded only needs a more careful definition of the domain $S_{(\bfx,\bfy)}$, replacing the parameter $h$ by a local value.
\end{proof}

\bigskip
The following bound provides an estimate of the translates of a discrete function $v$ as a function of $\normHb{v}$.

\begin{lemma}
Let $\theta_0 >0$ and let $\mesh$ be a triangulation of the computational domain $\Omega$ such that $\theta \geq \theta_0$, where $\theta$ is defined by \eqref{regul}; let $V_h$ be the space of Crouzeix-Raviart discrete functions associated to $\mesh$, as defined by \eqref{crdef}.
Let $v$ be a function of $V_h$; we denote by $\tilde v$ the extension by zero of $v$ to $\xR^d$.
Then the following estimate holds:
\[
\forall \bfeta \in \xR^d,\qquad \normLddR{\tilde v(\cdot + \bfeta) - \tilde v (\cdot)} \leq
c(\theta_0)\ |\bfeta|\ (|\bfeta| + h)\ \normHbd{v}.
\]
\label{trans_est}\end{lemma}
\begin{proof}
We follow the proof of a similar result for piecewise constant functions, namely \cite[Lemma 9.3, pp. 770-772]{eym-00-fin}.
Let $\bfeta \in \xR^d$ be given, $v$ be a Crouzeix-Raviart discrete function and $\tilde v$ its extension by zero to $\xR^d$.
With the definition of the function $\chi_\edge$ of Lemma \ref{ex_dedge}, the following identity holds for any $\bfx \in \xR^d$:
\[
\tilde v(\bfx+\bfeta)- \tilde v(\bfx) = 
\underbrace{\sum_{\edge \in \edges} \chi_\edge(\bfx,\bfx+\bfeta)\, [v](\bfy_{\bfx,\bfeta,\edge})}_{\displaystyle T_1(\bfx)}
+ \underbrace{\int_0^{1} \grad_h \tilde v(\bfx + s \bfeta) \cdot  \bfeta \ds}_{\displaystyle T_2(\bfx)}
\]
where $\bfy_{\bfx,\bfeta,\edge}$ stands for the intersection between the line issued from $\bfx$ and of direction $\bfeta$ and the hyperplane containing $\edge$.
Defining for each edge $\edge$ of the mesh a real positive number $d_\edge$ such that Lemma \ref{ex_dedge} holds, by the Cauchy-Schwarz inequality, we get for $T_1(\bfx)$:
\[
(T_1(\bfx))^2 \leq 
\left( \sum_{\edge \in \edges} \chi_\edge(\bfx,\bfx+\bfeta)\, \frac{[v](\bfy_{\bfx,\bfeta,\edge})^2}{d_\edge} \right)
\left( \sum_{\edge \in \edges} \chi_\edge(\bfx,\bfx+\bfeta)\  d_\edge \right)
\]
Integrating now over $\xR^d$, we thus obtain:
\[
\int_{\xR^d} (T_1(\bfx))^2 \dx \leq c_2(\theta_0)\ \left( |\bfeta| + h \right) 
\left( \sum_{\edge \in \edges} \int_{\xR^d} \chi_\edge(\bfx,\bfx+\bfeta)\, \frac{[v](\bfy_{\bfx,\bfeta,\edge})^2}{d_\edge} \dx \right)
\]
Let  $Q_{\edge,\eta} = \{ \bfx= \bfy + s \bfeta ; \bfy \in \edge \mbox{ and } s \in [-1,0]\}$. 
Noting that the function $\bfx \mapsto \chi_\edge(\bfx,\bfx+\bfeta)$ is in fact the characteristic function of $Q_{\edge,\eta}$, we get that:
\begin{eqnarray*}
\int_{\xR^d} \chi_\edge(\bfx,\bfx+\bfeta)\ \left(  [v](\bfy_{\bfx,\bfeta,\edge})\right)^2 \dx& =&\int_{Q_{\edge,\eta} }\left([v](\bfy_{\bfx,\bfeta,\edge})\right)^2 \dx\\
& = &\vert \bfn_\edge \cdot \bfeta \vert  \int_{-1}^0 \int_\edge\left(  [v](\bfy) \right)^2 \rm d \bfy  \ ds,
\end{eqnarray*}
where $\bfn_\edge$ is a unit normal vector to $\edge$. Therefore:
\[
\int_{\xR^d} (T_1(\bfx))^2 \dx \leq c_2(\theta_0)\ \left( |\bfeta| + h \right)\ |\bfeta| \, \sum_{\edge \in \edges} \frac {1} {d_\edge} \int_\edge \left(  [v](\bfy) \right)^2   \rm d \bfy,
\]
and thus, by choice of $d_\edge$:
\begin{equation}
\int_{\xR^d} (T_1(\bfx))^2 \dx \leq \frac{c_2(\theta_0)}{c_1(\theta_0)}\ \left( |\bfeta| + h \right)\ |\bfeta| \, \sum_{\edge \in \edges} \frac 1 {h_\edge} \int_\edge  \left(  [v](\bfy) \right)^2   \rm d \bfy.
\label{translates_1}\end{equation}
On the other hand, by the Cauchy-Schwarz inequality, we have for $T_2$:
\[
(T_2(\bfx))^2 \leq |\bfeta|^2 \int_0^{1} |\grad_h \tilde v(\bfx + s \bfeta)|^2 \ds,
\]
and thus, using the Fubini theorem and remarking that $\grad_h \tilde v$ vanishes outside $\Omega$, we  get:
\begin{equation}
\int_{\xR^d} (T_2(\bfx))^2 \dx \leq |\bfeta|^2 \normHbd{v}.
\label{translates_2}\end{equation}
The result then follows thanks to the inequality $| \tilde v(\bfx+\bfeta)- \tilde v(\bfx)|^2 \leq 2 (T_1(\bfx))^2 + 2 (T_2(\bfx))^2$,  to  the bounds \eqref{translates_1} and \eqref{translates_2} and to Lemma \ref{est_jumps}.
\end{proof}

We are now in position to state the following compactness result.
\begin{theorem}
Let $(v^{(m)})_{m \in \xN}$ be a sequence of functions satisfying the following assumptions:
\begin{enumerate}
\item $\forall m \in \xN$, there exists a triangulation of the domain $\mesh^{(m)}$ such that $v^{(m)} \in V_h^{(m)}$, where $V_h^{(m)}$ is the Crouzeix-Raviart space associated to $\mesh^{(m)}$ defined by \eqref{crdef}, and the parameter $\theta^{(m)}$ characterizing the regularity of $\mesh^{(m)}$ is bounded away from zero independently of $m$,
\item the sequence $(v^{(m)})_{m \in \xN}$ is uniformly bounded with respect to the broken Sobolev $\xHone$ semi-norm, \ie:
\[
\forall m \in \xN, \qquad \normHb{v^{(m)}} \leq C
\]
where the real number $C$ does not depend on $m$ and $\normHb{\cdot}$ stands for the broken Sobolev $\xHone$ semi-norm associated to $\mesh^{(m)}$ (with a slight abuse of notation, namely dropping for short the index $^{(m)}$ pointing the dependence of the norm with respect to the mesh).
\end{enumerate}
Then, possibly up to the extraction of a subsequence, the sequence $(v^{(m)})_{m \in \xN}$ converges strongly in $\xLtwo(\Omega)$ to a limit $\bar v$ such that $\bar v \in \xHone_0(\Omega)$.
\label{compactness}\end{theorem}
\begin{proof}
The result follows from the translates estimates of lemma \ref{trans_est}.
The compactness in $\xLtwo(\Omega)$ of the sequence is a consequence of the Kolmogorov theorem (see \eg\ \cite[theorem 14.1, p. 833]{eym-00-fin} for a statement of this result).
The fact that the limit belongs to $\xHone_0(\Omega)$ follows from the particular expression for the bound of the translates and is proven in \cite[theorem 14.2, pp. 833-834]{eym-00-fin}.
\end{proof}


\section{The numerical scheme} \label{sec:scheme}

Let $\rho^\ast$ be the mean density, \ie\ $\rho^\ast=M/|\Omega|$ where $|\Omega|$ stands for the measure of the domain $\Omega$.
We consider the following numerical scheme for the discretization of Problem (\ref{pbcont}):
\begin{subequations}
\begin{align}
& \forall \bfv \in \bfWh,  \,
\displaystyle\int_\Omega \grad_h \bfu : \grad_h \bfv \dx 
- \int_\Omega p\ \dive_h \bfv \dx = \int_\Omega \bff \cdot \bfv \dx, \label{pbdisc1}
\\&  
\forall K \in \mesh, \displaystyle \,
\sum_{\edge=K|L} \left( \vsp\, \varrho(p_K) - \vsm\ \varrho(p_L) \right) 
 + \underbrace{ h^\alpha\, |K|\, \left( \varrho(p_K)-\rho^\ast \right) }_{\displaystyle T_{{\rm stab},1}}    \label{pbdisc2}
 \\ & \hspace{10ex}+ \underbrace{\sum_{\edge=K|L} (h_K+h_L)^\beta\ \frac{|\edge|}{h_\edge}\ \left|\varrho(p_K) +\varrho(p_L) \right| 
 \ \left(\varrho(p_K) -\varrho(p_L) \right)}_{\displaystyle T_{{\rm stab},2}} =0, \nonumber
\end{align}\label{pbdisc}\end{subequations}
where $\vsp$ and $\vsm$ stands respectively for $\vsp = \max (\vs,\ 0)$ and $\vsm = -\min (\vs,\ 0)$ with $\vs=|\edge|\, \bfu_\edge \cdot \bfn_{KL}=\vsp- \vsm$.
Note that the upwinded convection term $\sum_{\edge=K|L} \left( \vsp\, \varrho(p_K) - \vsm\ \varrho(p_L) \right) $ may also be written:   $\sum_{\edge=K|L} \vs \rho_\sigma$, with
\begin{equation}
 \rho_\sigma = \begin{cases} \rho_K \text{ if } \vs \ge 0,\\ \rho_L \text{ otherwise.} \end{cases}
\label{upwindchoice}
\end{equation}
Equation \eqref{pbdisc1} may be considered as the standard finite element discretization of \eqref{pbcont1}.
Since the pressure is piecewise constant, the finite element discretization of \eqref{pbcont2}, \ie\ the mass balance, is similar to a finite volume formulation, in which we introduce the standard first-order upwinding and two stabilizing terms.
The first one, \ie\ $T_{{\rm stab},1}$, guarantees that the integral of the density over the computational domain is always $M$ (this can easily be seen by summing the second relation for $K\in\mesh$).
The second one, \ie \ $T_{{\rm stab},2}$, is useful in the convergence analysis.
It may be seen as a finite volume analogue of a continuous term of the form $\dive \left( |\rho| \grad \rho\right) $ weighted by a mesh-dependent coefficient tending to zero as $h^\beta$; note, however, that $h_\edge$ is not the distance which is  usually encountered in the finite volume discretization of diffusion terms; consequently, the usual restriction for the mesh when diffusive terms are to be approximated by the two-points finite volume method (namely,the Delaunay condition) is not required here.
We suppose that $\alpha \geq 1$ and the convergence analysis uses $0 < \beta <2$.

\begin{remark}
At first glance, leaving the weight $|\rho|$ out, the stabilization term $T_{{\rm stab},2}$ may look as a Brezzi-Pitk\"{a}ranta regularisation \cite{bre-84-sta}, as used in \cite{eym-06-sta} for stabilizing the colocated approximation of the Stokes problem, which would be rather puzzling since we use here an {\em inf-sup} stable pair of approximation spaces.
However, using the equation of state \eqref{eos}, we obtain:
\[
T_{{\rm stab},2}= \machd^2
\sum_{\edge=K|L} (h_K+h_L)^\beta\ \frac{|\edge|}{h_\edge}\ \left|p_K + p_L \right| 
\ \left( p_K -p_L \right)
\]
which shows, since $\machd^2=\mach^4$, that this term rapidly vanishes when approaching the incompressible limit.
\end{remark}


\section{Existence of a solution and a priori estimates}\label{sec:apriori}

The existence of a solution to \eqref{pbdisc} follows, with minor changes to cope with the diffusion stabilization term, from the theory developed in \cite[section 2]{gal-08-unc}.
In this latter paper, it is obtained for fairly general equations of state by a topological degree argument.
We only give here the obtained result, together with a proof of the {\it a priori} estimates verified by the solution, and we refer to \cite{gal-08-unc} for the proof of existence.

\begin{theorem}
Let $\theta_0 >0$ and let $\mesh$ be a triangulation of the computational domain $\Omega$ such that $\theta \geq \theta_0$, where $\theta$ is defined by \eqref{regul}.
Problem \eqref{pbdisc} admits at least one solution $(\bfu,p) \in \bfWh \times L_h$;  any possible solution satisfies $p_K >0$, $\forall K \in \mesh$ and:
\begin{equation}
\normHb{\bfu} + \normLd{p} + \normLd{\rho} + \snormdisc{\rho} \leq C(\bff,M)
\label{apriori}\end{equation}
where $C(\bff,M) \in \xR$ only depends on $\Omega$, $\machd$, $\bff$, $M$ and $\theta_0$.
\label{prop:apriori}\end{theorem}

\begin{proof}
Let $(\bfu,p) \in \bfWh \times L_h$ be a solution to \eqref{pbdisc}.
Let $\rho_K = \varrho(p_K)$ for any $K \in \mesh$, and let $\rho$ denote the vector $(\rho_K)_{K \in \mesh}$. 
A natural ordering of the equations and unknowns in  (\ref{pbdisc2}) leads to a linear system of the form $M \rho = c$, where $c \in \xR^N$,  $N$ is the number discretization   cells,  $c \in \xR^N$, $c >0$, and where $M$ is an $M$--matrix (in particular $M^{-1} \ge 0$ and  $M^{-t} \ge 0$). 
Therefore the $i$-th component of $\rho$ reads $\rho_i = M^{-1}c\cdot e_i = c\cdot M^{-t} e_i$ where $e_i$ is the  $i$-th canonical basis vector of  $\xR^N$.
Since   $M^{-t} \ge 0$, we get  $M^{-t} e_i \ge 0$, and since $M^{-t} e_i \not = 0$, this proves that $\rho_i >0$, which, in turns,  yields $p_K >0$, $\forall K \in \mesh$. 
Let us then prove the estimate \eqref{apriori}.
To this end, we take $\bfv = \bfu$ in   \eqref{pbdisc1} and obtain:
\begin{equation}
\int_\Omega |\grad_h \bfu|^2 \dx 
- \int_\Omega p\ \dive_h \bfu \dx = \int_\Omega \bff \cdot \bfu \dx.
\label{apriori_1}\end{equation}
Let us then multiply \eqref{pbdisc2} by $\machd^{-1} [1 + \log (\rho_K)]$ (see remark \ref{logrho} below for an explanation of this choice) and we sum over $K \in \mesh$; dropping the terms which vanish by conservativity, we then obtain:
\[
\begin{array}{l}
T_1 + T_2 + T_3 =0 \qquad \mbox{with:}
\\[2ex] \hspace{3ex} \left| \begin{array}{l}\displaystyle
T_1= \machd^{-1} \sum_{K\in\mesh} \log (\rho_K) \sum_{\edge=K|L} \left( \vsp\, \rho_K - \vsm\ \rho_L\right),
\\[3ex] \displaystyle
T_2= \machd^{-1} \, h^\alpha \sum_{K\in\mesh} |K|\, \left[ 1 +\log (\rho_K)\right]\, \left[ \rho_K-\rho^\ast \right],
\\[3ex] \displaystyle
T_3=\machd^{-1} \sum_{K\in\mesh} \log ( \rho_K) 
\sum_{\edge=K|L} (h_K+h_L)^\beta\ \frac{|\edge|}{h_\edge} \left( \rho_K + \rho_L \right) \ \left( \rho_K - \rho_L \right),
\end{array} \right.
\end{array}
\]
where  the term $|\varrho(p_K) +\varrho(p_L)|$ has been replaced by $\varrho(p_K) +\varrho(p_L)$ in \eqref{pbdisc2}, thanks to the positivity of the pressure.
Let us first write $T_1$ as:
\[
T_1=\machd^{-1} \sum_{K\in\mesh} \log (\rho_K) \sum_{\edge=K|L} \vs\, \rho_\edge,
\]
where  $\rho_\edge$ is the upwind choice defined by \eqref{upwindchoice}.
Adding and substracting the same quantity, $T_1$ equivalently reads:
\[
\begin{array}{l} \displaystyle
T_1
=\machd^{-1} \sum_{K\in\mesh} \rho_K \sum_{\edge=K|L} \vs
+\machd^{-1}  \sum_{K\in\mesh} \sum_{\edge=K|L} \vs (\rho_\edge \log (\rho_K) -\rho_K).
\end{array}
\]
In the first summation, we recognize $\displaystyle \int_\Omega p\ \dive_h \bfu \dx$.
A reordering of the second summation yields:
\[
T_1=\int_\Omega p\ \dive_h \bfu \dx + \machd^{-1} \sum_{\stackrel{\scriptstyle \edge \in \edgesint,}{\scriptstyle \edge=K|L}} \vs \left[(\rho_\edge \log(\rho_K) - \rho_K) -(\rho_\edge \log(\rho_L) - \rho_L)\right].
\]
Let $\bar \rho_\edge$ be defined by  $\begin{cases} \bar \rho_\edge =\rho_K=\rho_L \text{ if } \rho_K=\rho_L, \\\bar \rho_\edge \log(\rho_K) - \rho_K = \bar \rho_\edge \log(\rho_L) - \rho_L \text{ otherwise}.
\end{cases}$ \\
With this notation, we get:
\[
T_1=\int_\Omega p\ \dive_h \bfu \dx +
\machd^{-1}  \sum_{\stackrel{\scriptstyle \edge \in \edgesint,}{\scriptstyle \edge=K|L}} \vs 
\ (\rho_\edge - \bar \rho_\edge)\ (\log(\rho_K) - \log(\rho_L)).
\]
In the last summation, we can, without loss of generality, choose the orientation of each edge in such a way that $\vs \geq 0$.
With this convention, the term in the summation reads $\vs \ (\rho_K - \bar \rho_\edge)\ (\log(\rho_K) - \log(\rho_L))$, and is non-negative thanks to the fact that $\rho_\edge \in [\min (\rho_K,\rho_L),\ \max (\rho_K,\rho_L)]$ and the $\log$ function  is increasing.
We thus finally obtain:
\begin{equation}
T_1 \geq \int_\Omega p\ \dive_h \bfu \dx.
\label{apriori_2}\end{equation}
Let us now turn to the estimate of $T_2$.
Since the function $z \mapsto z \log(z)$ is convex for positive $z$ and its derivative is $z \mapsto 1 + \log(z)$, we simply have:
\begin{equation}
T_2 \geq \machd^{-1} \, h^\alpha \sum_{K\in\mesh} |K|\, \left[ \rho_K \log(\rho_K)-\rho^\ast \log(\rho^\ast)\right].
\label{apriori_3}\end{equation}
Finally, reordering the sums, the term $T_3$ reads:
\[
T_3= \machd^{-1} \sum_{\stackrel{\scriptstyle \edge \in \edgesint,}{\scriptstyle \edge=K|L}}
(h_K+h_L)^\beta\ \frac{|\edge|}{h_\edge}\ \left(\rho_K +\rho_L \right) \ \left(\rho_K -\rho_L \right)
\ \left(\log(\rho_K) -\log(\rho_L) \right).
\]
By concavity of the $\log$ function, we have:
\[
|\log(\rho_K) -\log(\rho_L)| \geq \frac 1 {\max (\rho_K,\rho_L)} \ |\rho_K -\rho_L|,
\]
and thus:
\begin{equation}
T_3 \geq \machd^{-1} \sum_{\stackrel{\scriptstyle \edge \in \edgesint,}{\scriptstyle \edge=K|L}}
(h_K+h_L)^\beta\ \frac{|\edge|}{h_\edge}\ \left(\rho_K -\rho_L \right)^2.
\label{apriori_4}\end{equation}
Summing equations \eqref{apriori_1}--\eqref{apriori_4} and using Young's inequality, we obtain:
\[
\normHb{\bfu} + \machd^{-1/2} \snormdisc{\rho} \leq C(\bff,M).
\]
Furthermore, summing  \eqref{pbdisc2} over $K \in \mesh$, we obtain that the mean value of the pressure $p_{\rm m}$ is given by:
\[
p_{\rm m}=\frac 1 {|\Omega|} \int_\Omega p \dx = \machd^{-1} \, \rho^\ast.
\]
Using the {\it inf-sup} stability of the discretization, we get on the other hand:
\[
\begin{array}{ll} \displaystyle
\normLd{p-p_{\rm m}} 
& \displaystyle
\leq \frac 1 {c_{\rm i}}
\sup_{\bfv \in \bfWh} \frac 1 {\normHb{\bfv}} \int_\Omega p \ \dive_h \bfv \dx 
\\[3ex] & \displaystyle
= \frac 1 {c_{\rm i}}\sup_{\bfv \in \bfWh} \frac 1 {\normHb{\bfv}} \int_\Omega \left( \grad_h \bfu : \grad_h \bfv - \bff \cdot \bfv \right) \dx,
\end{array}
\]
and the control of $\normLd{p}$ (or, equivalently, $\normLd{\rho}$) follows from the estimate for $\normHb{\bfu}$.
\end{proof}

\begin{remark}[On the choice of $(\log(\rho_K))_{K \in \mesh}$ as test function]
At first glance, the choice of $\log(\rho_K)$ to multiply \eqref{pbdisc2} in the preceding proof may seem rather puzzling.
In fact, this computation is a particular case of the so-called "elastic potential identity", which is well-known in the continuous setting and is central in {\it a priori} estimates for the compressible Navier-Stokes equations \cite{pll-98-mat,nov-04-int,fei-04-dyn}.
An analogous identity is proven at the discrete level, for the same discretization as here, in \cite[theorem 2.1]{gal-08-unc}.

For the particular case under consideration, an elementary explanation of this choice may be given.
Indeed, it is crucial in the above proof that the quantity $\bar \rho_\edge$ lies in the interval $[\min (\rho_K,\rho_L),\ \max (\rho_K,\rho_L)]$.
Let us suppose, without loss of generality, that $0< \rho_K < \rho_L$ and that, instead of the $\log$  function, the computation is performed with a non-specified increasing ond continuously differentiable function $f$; then we get:
\[
\bar \rho_\edge = \frac{\rho_L-\rho_K}{f(\rho_L)-f(\rho_K)}.
\]
The condition $\bar \rho_\edge \geq \rho_K$ is equivalent to:
\[
\frac 1 {\rho_K} \geq \frac{f(\rho_L)-f(\rho_K)}{\rho_L-\rho_K},
\]
which is verified for $f=\log$ by concavity of the latter and, letting $\rho_L$ tend to $\rho_K$, yields $f'(x) \leq 1/x$.
Conversely, the condition $\bar \rho_\edge \leq \rho_L$ yields:
\[
\frac 1 {\rho_L} \leq \frac{f(\rho_L)-f(\rho_K)}{\rho_L-\rho_K},
\]
which, once again, is verified by the $\log$ function, and now implies $f'(x) \geq 1/x$.

This limitation for the choice of the test function is the reason for the expression of the stabilizing diffusion  term. 
\label{logrho}\end{remark}


\section{Convergence analysis}\label{sec:conv}

In this section, we prove the following convergence result.

\begin{theorem}
Let $(\mesh^{(m)})_{m \in \xN}$ be a sequence of triangulations of $\Omega$ such that $h^{(m)}$ tends to zero when $m$ tends to $+ \infty$.
Let us assume that this sequence is regular in the sense that there exists $\theta_0 >0$ such that $\theta^{(m)} \geq \theta_0,\ \forall m \in \xN$, where $\theta^{(m)}$ is defined by \eqref{regul}.
For $m \in \xN$, we denote by $\bfWh^{(m)}$ and $L_h^{(m)}$  the discrete  velocity and pressure spaces  associated to $\mesh^{(m)}$ and by $(\bfu^{(m)},p^{(m)}) \in \bfWh^{(m)}\times L_h^{(m)}$ the corresponding solution to \eqref{pbdisc}, with $\alpha\geq 1$ and $0<\beta<2$.
Then, up to a subsequence, the sequence $(\bfu^{(m)})_{m \in \xN}$ strongly converges to a limit $\bar \bfu$ in $\xLtwo(\Omega)^d$ and $(p^{(m)})_{m \in \xN}$ converges to $\bar p$ weakly in $\xLtwo(\Omega)$, where the pair $(\bar \bfu,\bar p)$ is a solution to \eqref{pbcont} in the following weak sense:
\begin{subequations} \begin{align} & \displaystyle \nonumber
\bar \bfu \in \xHone_0(\Omega)^d,\ \bar p \in \xLtwo(\Omega) \mbox{ and }:
\\ & \displaystyle \qquad \label{conv_1}
\int_\Omega \grad \bar \bfu : \grad \bfpsi \dx - \int_\Omega \bar p \ \dive \bfpsi \dx = \int_\Omega \bff \cdot \bfpsi \dx,
\qquad &
\forall \bfpsi \in \xCinfty_c(\Omega)^d,
\\ & \displaystyle \qquad \label{conv_2}
\int_\Omega \bar p \, \bar \bfu \cdot \grad \psi \dx =0,
&
\forall \psi \in \xCinfty_c(\Omega),
\\ & \displaystyle \qquad \label{conv_3}
\int_\Omega \varrho(\bar p)=M.
\end{align} \label{conv}\end{subequations}
\end{theorem}
\begin{proof}
The proof is divided in three steps: we first show the existence of the limits $\bar \bfu$ and $\bar p$, then we pass to the limit in \eqref{pbdisc1} and \eqref{pbdisc2}.
Since the equation of state is linear, the last equation is then a straightforward consequence of the weak convergence in $\xLtwo(\Omega)$ of the (sub)sequence $(p^{(m)})_{m \in \xN}$ to $\bar p$.

\vspace{2ex} \noindent \underline{\it Step 1: existence of a limit.}
 
\medskip
\noindent By the {\it a priori} estimates of theorem \ref{prop:apriori}, we know that:
$
\forall m \in \xN, \normHb{\bfu^{(m)}} \leq C(\bff,M).
$
The compactness in $\xLtwo(\Omega)^d$ of the sequence $(\bfu^{(m)})_{m\in\xN}$, together with the fact that the limit $\bar \bfu$ lies in $\xHone_0(\Omega)^d$, thus follows by applying theorem \ref{compactness} to each component $u^{(m)}_i,\ 1\leq i \leq d$.
Once again by theorem \ref{prop:apriori}, we have:
$
\forall m \in \xN,\quad \normLd{p^{(m)}} \leq C(\bff,M).
$
which is sufficient to ensure a weak convergence in $\xLtwo(\Omega)$ of the sequence $(p^{(m)})_{m\in\xN}$ to $\bar p \in \xLtwo(\Omega)$.

\vspace{2ex} \noindent \underline{\it Step 2: passing to the limit in \eqref{pbdisc1}.}

\medskip
\noindent Let $\bfpsi$ be a function of $\xCinfty_c(\Omega)^d$.
We denote by $\bfpsi^{(m)}$ the interpolation of $\bfpsi$ in $\bfWh^{(m)}$, \ie\ $\bfpsi^{(m)}=r_h^{(m)} \bfpsi$ where the operator $r_h^{(m)}$ is defined by \eqref{def_rh}.
Taking $v =\bfpsi^{(m)}$ in \eqref{pbdisc1}, we get:
\[
\int_\Omega \grad_h \bfu^{(m)} : \grad_h \bfpsi^{(m)} \dx - \int_\Omega p^{(m)} \ \dive_h \bfpsi^{(m)} \dx = 
\int_\Omega \bff \cdot \bfpsi^{(m)} \dx, \, \forall m \in \xN.
\]
Since the considered interpolation operator preserves the divergence \eqref{interp_1}, we have:
\[
\int_\Omega p^{(m)} \ \dive_h \bfpsi^{(m)} \dx =\int_\Omega p^{(m)} \ \dive \bfpsi \dx
\longrightarrow \int_\Omega \bar p \ \dive \bfpsi \dx \quad \mbox{ as } m \longrightarrow +\infty.
\]
By the approximation properties of the interpolation operator \eqref{interp_3} invoked component by component, we have:
\[
\int_\Omega \bff \cdot \bfpsi^{(m)} \dx \longrightarrow \int_\Omega \bff \cdot \bfpsi \dx \quad \mbox{ as } m \longrightarrow \infty.
\]
Finally, we have:
\[
\begin{array}{l} \displaystyle
\int_\Omega \grad_h \bfu^{(m)} : \grad_h \bfpsi^{(m)} \dx = 
\\ \displaystyle \hspace{20ex}
\underbrace{\int_\Omega \grad_h \bfu^{(m)} : \grad_h (\bfpsi^{(m)}-\bfpsi) \dx}_{\displaystyle T_1^{(m)}}
+ \underbrace{\int_\Omega \grad_h \bfu^{(m)} : \grad \bfpsi \dx}_{\displaystyle T_2^{(m)}}
\end{array}
\]
Once again by \eqref{interp_3}, the term $T_1^{(m)}$ obeys the following estimate:
\[
|T_1^{(m)}| \leq \normHb{\bfu^{(m)}} \normHb{\bfpsi^{(m)}-\bfpsi} \leq c(\theta_0) \, h^{(m)} \normHb{\bfu^{(m)}}
\snormHdeuxdom{\bfpsi}{\Omega},
\]
and thus tends to zero as $m$ tends to $+ \infty$.
Integrating by parts over each control volume, the term $T_2^{(m)}$ reads:
\[
T_2^{(m)}= -\int_\Omega \bfu^{(m)} \cdot \Delta \bfpsi \dx 
+ \sum_{\edge \in \edgesint^{(m)}} \int_\edge [\bfu^{(m)}]\, \grad \bfpsi \cdot \bfn_\edge \dedge,
\]
where $\bfn_\edge$ is a normal vector to $\edge$, with the same orientation as that of the jump through $\edge$.
Applying Lemma \ref{jump_reg} for each component of $\grad \bfpsi$, $a_\edge$ being the relevant component of the normal vector $n_\edge$, we get:
\[
\sum_{\edge \in \edgesint^{(m)}} \int_\edge [\bfu^{(m)}]\, \grad \bfpsi \cdot \bfn_\edge \dedge \leq 
c(\theta_0) \, h^{(m)} \normHb{\bfu^{(m)}} \snormHdeuxdom{\bfpsi}{\Omega},
\]
and thus tends to zero, while the first one tends to $-\int_\Omega \bar \bfu \cdot \Delta \bfpsi \dx$ as $m$ tends to $+ \infty$.
Since $\bar \bfu \in \xHone_0(\Omega)^d$, we may integrate by parts, and collecting the limits, we obtain \eqref{conv_1}.

\vspace{2ex} \noindent \underline{\it Step 3: passing to the limit in \eqref{pbdisc2}.}

\medskip
\noindent Let $\psi$ be a function of $\xCinfty_c(\Omega)$.
Multiplying the second equation of \eqref{pbdisc} by $ 1/ {|K|}\ \int_K \psi(\bfx)\dx$ and summing over $K\in\mesh$ yields $T_3^{(m)} + T_4^{(m)} + T_5^{(m)} =0, \ \forall m \in \xN$, with:
\[
\begin{array}{l} \displaystyle
T_3^{(m)}= \sum_{K \in \mesh^{(m)}} \frac 1 {|K|} \left( \sum_{\edge=K|L} \vs^{(m)}\  \rho_\edge^{(m)} \right)\ \int_K \psi \dx
\\[5ex] \displaystyle
T_4^{(m)}= (h^{(m)})^\alpha \sum_{K \in \mesh^{(m)}}|K| \left( \rho_K^{(m)}-\rho^\ast \right)\ \psi_K 
\\[4ex] \displaystyle
T_5^{(m)}= \sum_{K \in \mesh^{(m)}} 
\left(\sum_{\edge=K|L} (h_K+h_L)^\beta\ \frac{|\edge|}{h_\edge}\ \left(\rho_K^{(m)} +\rho_L^{(m)} \right)
\ \left(\rho_K^{(m)} -\rho_L^{(m)} \right)\right) \ \psi_K,
\end{array}
\]
where $\rho_\edge^{(m)}$ is defined by \eqref{upwindchoice} and $\psi_K$ stands for the mean value of $\psi$ over $K$.
Let $\bfq^{(m)} \in \bfWh$ be defined as $\bfq^{(m)}(\bfx)=\sum_{\edge \in \edgesint^{(m)}} \bfu_\edge^{(m)}\, \rho_\edge^{(m)}\, \crou_\edge(\bfx)$, where $\crou_\edge$ is the Crouzeix-Raviart basis function associated to $\edge$.
The divergence of $\bfq^{(m)}$ is a piecewise constant function and reads:
\[
\forall K \in \mesh^{(m)}, \qquad
\dive \bfq^{(m)}=\frac 1 {|K|} \sum_{\edge=K|L} \vs^{(m)}\  \rho_\edge^{(m)} \qquad \mbox{a.e. in } K.
\]
We thus have for $T_3^{(m)}$:
\[
T_3^{(m)}=\sum_{K \in \mesh^{(m)}} \int_K \psi\ \dive \bfq^{(m)} \dx.
\]
Integrating by parts over each control volume, we get:
\[
\begin{array}{ll} \displaystyle
T_3^{(m)}
& \displaystyle
= -\int_\Omega \grad \psi \cdot \bfq^{(m)} \dx +
\sum_{\edge \in \edgesint^{(m)}} \int_\edge \psi \ [\bfq^{(m)}] \cdot \bfn_\edge \dedge
\\ & \displaystyle
= -\int_\Omega \grad \psi \cdot (\rho^{(m)}\,\bfu^{(m)}) \dx 
\\ & \displaystyle \hspace{10ex}
+ \underbrace{\sum_{\edge \in \edgesint^{(m)}} \int_\edge \psi \ [\bfq^{(m)}] \cdot \bfn_\edge \dedge}_{\displaystyle T_6^{(m)}}
+ \underbrace{\int_\Omega \grad \psi \cdot (\bfq^{(m)}-\rho^{(m)}\,\bfu^{(m)}) \dx}_{\displaystyle T_7^{(m)}}.
\end{array}
\]
By the respectively weak and strong convergence of $(\rho^{(m)})_{m\in\xN}$ and $(\bfu^{(m)})_{m \in \xN}$ to $\bar \rho$ and $\bar \bfu$ in $\xLtwo(\Omega)$ and $\xLtwo(\Omega)^d$, we have:
\[
\int_\Omega \grad \psi \cdot (\rho^{(m)}\,\bfu^{(m)}) \dx \longrightarrow \int_\Omega \grad \psi \cdot (\bar \rho\,\bar \bfu ) \dx
\qquad \mbox{ as } m \longrightarrow + \infty.
\]
By the definition of $\bfq^{(m)}$, the term $T_6^{(m)}$ reads:
\[
T_6^{(m)} 
= \sum_{\edge \in \edgesint^{(m)}} \int_\edge \psi 
\ [\sum_{\edge' \in \edgesint^{(m)}} \bfu_{\edge'}^{(m)}\, \rho_{\edge'}^{(m)}\, \crou_{\edge'}(\bfx)] \cdot \bfn_\edge \dedge
=T_8^{(m)} + T_9^{(m)}.
\]
with:
\[
\begin{array}{ll} \displaystyle
T_8^{(m)} 
& \displaystyle =
\sum_{\edge \in \edgesint^{(m)}} \int_\edge \psi 
\ \rho_\edge^{(m)}[\sum_{\edge' \in \edgesint^{(m)}} \bfu_{\edge'}^{(m)}\,  \crou_{\edge'}(\bfx)] \cdot \bfn_\edge \dedge
\\ & \displaystyle =
\sum_{\edge \in \edgesint^{(m)}} \rho_\edge^{(m)} \int_\edge \psi \ [\bfu^{(m)}] \cdot \bfn_\edge \dedge
\\ \displaystyle
T_9^{(m)}
& \displaystyle =
\sum_{\edge \in \edgesint^{(m)}} \int_\edge \psi 
\ [\sum_{\edge' \in \edgesint^{(m)}\setminus \{\edge\}} \bfu_{\edge'}^{(m)}\, (\rho_{\edge'}^{(m)}-\rho_\edge^{(m)})\, \crou_{\edge'}(\bfx)] \cdot \bfn_\edge \dedge.
\end{array}
\]
Since the integral of the jump of a Crouzeix-Raviart function over an internal edge of the mesh vanishes, the term $T_8^{(m)}$ can be estimated as follows:
\[
|T_8^{(m)}| \leq c_\psi\, h^{(m)} \sum_{\edge \in \edgesint^{(m)}} \rho_\edge^{(m)} 
\int_\edge \left| [\bfu^{(m)}] \cdot \bfn_\edge \right| \dedge,
\]
where $c_\psi$ only depends on $\psi$.
Using the Cauchy-Schwarz inequality then yields:
\[
\begin{array}{ll}\displaystyle
|T_8^{(m)}|
& \displaystyle
\leq c_\psi\, h^{(m)} \sum_{\edge \in \edgesint^{(m)}} |\edge|^{1/2}\ \rho_\edge^{(m)} 
\Bigl( \int_\edge | \ [\bfu^{(m)}  ] \ |^2 \dedge \Bigr)^{1/2}
\\ & \displaystyle
\leq c_\psi\, h^{(m)}  \Bigl( \sum_{\edge \in \edgesint^{(m)}} h_\edge\, |\edge| \, (\rho_\edge^{(m)})^2 \Bigr)^{1/2}
\Bigl(\sum_{\edge \in \edgesint^{(m)}} \frac 1 {h_\edge} \int_\edge |\ [\bfu^{(m)}] \ |^2 \dedge \Bigr)^{1/2}.
\end{array}
\]
By the regularity of the mesh, the first summation is bounded by $\normLd{\rho^{(m)}}$ while, by Lemma \ref{est_jumps}, the second one is bounded by $c(\theta_0) \normHbd{\bfu^{(m)}}$.
Let us now turn to the study of $T_9^{(m)}$.
Since, for $\edge' \in \edgesint^{(m)}\setminus \{\edge\}$, the integral of $\crou_{\edge'}$ over $\edge$ vanishes, and since the functions $\crou_\edge$ are bounded (namely $\vert \crou_\edge \vert \le 1$ in 2D,  $\vert \crou_\edge \vert \le 2$ in 3D) we get:
\[
\int_\edge \psi\ (\rho_{\edge'}^{(m)}-\rho_\edge^{(m)})\, [\crou_{\edge'}(\bfx)]\ \bfu_{\edge'}^{(m)}\cdot \bfn_\edge \dedge
\leq c_\psi\, h_\edge\ |\edge|\ |\rho_{\edge'}^{(m)}-\rho_\edge^{(m)}|\ |\bfu_{\edge'}^{(m)}|,
\]
where $c_\psi$ still only depends on $\psi$.
Since the function $\crou_{\edge'}$ is non-zero over $\edge=K|L$ only when $\edge'$ belongs to the edges of $K$ or $L$, only a limited number of terms are non-zero in $T_9^{(m)}$, in such a way that the difference $\rho_{\edge'}^{(m)}-\rho_\edge^{(m)}$ only involves two neigbouring cells or two cells sharing the same neighbour.
Splitting the difference in this last case, using the previous inequality and the regularity of the mesh (in particular the fact that the ratio of the size of two neighbouring cells is bounded) and reordering the sums, we get for $T_9^{(m)}$ an estimate of the form:
\[
|T_9^{(m)}| \leq c \sum_{\edge \in \edgesint^{(m)}} h_\edge^d \ |\bfu_\edge^{(m)}|
\sum_{\stackrel{\scriptstyle \edge' \in {\mathcal N}_\edge}{\scriptstyle (\edge'=K|L)}} |\rho_K^{(m)}-\rho_L^{(m)}|,
\]
where the positive real number $c$ only depends on $\psi$ and the regularity of the mesh and, thanks to this regularity, the set ${\mathcal N}_\edge$ is such that a given edge $K|L$ only appears in this sum a number of times bounded independently of $m$.
Thus, thanks to the Cauchy-Schwarz inequality, we have:
\[
|T_9^{(m)}|^2 \leq
c\  \Bigl(\sum_{\edge \in \edgesint^{(m)}} h_\edge^d\ |\bfu_\edge^{(m)}|^2 \Bigr)
\ \Bigl( \sum_{\stackrel{\scriptstyle \edge \in \edgesint^{(m)}}{\scriptstyle (\edge=K|L)}} h_\edge^d
\ \left(\rho_K^{(m)}-\rho_L^{(m)}\right)^2 \Bigr).
\]
By the regularity of the mesh, the first term of this product is controlled by $\normLd{\bfu^{(m)}}$ and the second one by $(h^{(m)})^{2-\beta} \snormdisc{\rho^{(m)}}$.
Consequently, thanks to estimate \eqref{apriori}, both $T_8^{(m)}$ and $T_9^{(m)}$ and thus also $T_6^{(m)}$ tend to zero as $m$ tends to $+\infty$, for any $\beta<2$.\\[2ex]
Let us then examine the term $T_7^{(m)}$:
\[
T_7^{(m)}=\sum_{K\in\mesh^{(m)}} \int_K \sum_{\edge=K|L} (\rho^{(m)}_\edge - \rho^{(m)}_K) \, \crou_\edge(\bfx) \, \bfu^{(m)}_\edge \cdot \grad \psi(\bfx) \dx
\]
Since $\grad \psi$ is bounded in $\xLinfty(\Omega)^d$, and since the functions $\crou_\edge$ are bounded, we get:
\[
|T_7^{(m)}|
\leq c_\psi \sum_{K\in\mesh^{(m)}} |K| \ \sum_{\edge=K|L} |\rho^{(m)}_\edge - \rho^{(m)}_K|\  |\bfu^{(m)}_\edge|.
\]
Reordering the summations and using the Cauchy-Schwarz inequality yields:
\[
|T_7^{(m)}| 
\leq c_\psi \sum_{\stackrel{\scriptstyle \edge \in \edgesint^{(m)}}{\scriptstyle (\edge=K|L)}}
(|K|+|L|)\ |\rho^{(m)}_K - \rho^{(m)}_L|\  |\bfu^{(m)}_\edge| \leq c_\psi\ \left(T_{10}^{(m)}\right)^{1/2} \left(T_{11}^{(m)}\right)^{1/2},
\]
with:
\[
\begin{array}{ll}\displaystyle
T_{10}^{(m)}= \sum_{\stackrel{\scriptstyle \edge \in \edgesint^{(m)}}{\scriptstyle (\edge=K|L)}} (|K|+|L|)\ |\bfu^{(m)}_\edge|^2
\\ \displaystyle
T_{11}^{(m)}=\sum_{\stackrel{\scriptstyle \edge \in \edgesint^{(m)}}{\scriptstyle (\edge=K|L)}}
h_\edge\, (h_K+h_L)^{(1-\beta)}\ \frac{|K|+|L|}{|\edge|\,(h_K+h_L)}
\ (h_K+h_L)^\beta\ \frac{|\edge|}{h_\edge}\ (\rho^{(m)}_K - \rho^{(m)}_L)^2.
\end{array}
\]
Once again reordering the summation, we get:
\[
T_{10}^{(m)}= \sum_{K \in \mesh}  |K|\ \sum_{\sigma \in \edges_K}|\bfu^{(m)}_\edge|^2,
\]
and thus, the term $T_{10}^{(m)}$ is controlled by $\normLd{\bfu^{(m)}}$, and $T_{11}^{(m)}$ is controlled by $(h^{(m)})^{2-\beta} \snormdisc{\rho^{(m)}}$.
By the {\it a priori} estimate \eqref{apriori}, $T_7^{(m)}$ thus tends to zero for any $\beta <2$.\\[2ex]
We now turn to the terms $T_4^{(m)}$ and $T_5^{(m)}$.
Since the sequence $(\rho^{(m)})_{m \in \xN}$ is bounded in $\xLtwo(\Omega)$, the term $T_4^{(m)}$ tends to zero for any $\alpha>0$.
Reordering the summation in $T_5^{(m)}$, we get:
\[
T_5^{(m)}=\sum_{\stackrel{\scriptstyle \edge \in \edgesint^{(m)}}{\scriptstyle (\edge=K|L)}}
(h_K+h_L)^\beta\ \frac{|\edge|}{h_\edge}\ \left(\rho_K^{(m)} +\rho_L^{(m)} \right)\ \left(\rho_K^{(m)}-\rho_L^{(m)}\right)
\ \left(\psi_K -\psi_L \right).
\]
By regularity of $\psi$, $|\psi_K -\psi_L| \leq c_\psi\ (h_K+h_L)$ and thus:
\[
|T_5^{(m)}| \leq c_\psi \sum_{\stackrel{\scriptstyle \edge \in \edgesint^{(m)}}{\scriptstyle (\edge=K|L)}}
(h_K+h_L)^{\beta+1}\ \frac{|\edge|}{h_\edge}\ \left(\rho_K^{(m)} +\rho_L^{(m)} \right)\ \left|\rho_K^{(m)}-\rho_L^{(m)}\right|.
\]
Using the Cauchy-Schwarz inequality, we obtain:
\[
|T_5^{(m)}| \leq c_\psi h^{\beta/2}
\left( \sum_{\stackrel{\scriptstyle \edge \in \edgesint^{(m)}}{\scriptstyle (\edge=K|L)}}
(h_K+h_L)^{2}\ \frac{|\edge|}{h_\edge}\ \left(\rho_K^{(m)} +\rho_L^{(m)} \right)^2 \right)^{1/2}
\ \snormdisc{\rho^{(m)}},
\]
which, once again since the sequence $(\rho^{(m)})_{m \in \xN}$ is bounded in $\xLtwo(\Omega)$, tends to zero by regularity of the mesh for $\beta > 0$.
The proof is thus complete.
\end{proof}


\section{Discussion}

To our knowledge, the convergence analysis performed in this paper seems to be the first result of this kind for the compressible Stokes problem (and, of course, more widely, for the compressible Navier-Stokes equations).
Beside the convergence of the scheme, it also provides an existence result for solutions of the continuous problem, which could also be derived from the continuous existence theory  ingredients for the steady Navier-Stokes equations (see \cite{nov-04-int} and references therein), but does not seem to be a direct consequence of the published literature: existence of strong solutions of the Navier-Stokes equations is known only for small data (\eg\ \cite{val-87-exi}) and existence of weak solutions is only proven for a particular class of equations of state (typically, $p=\rho^\gamma$ with $\gamma >3/2$), this limitation being due to the presence of the convection term.

\medskip
A puzzling fact is that the present theory relies on two ingredients which are usually not present in actual implementations.
Firstly, the stabilisation term $T_{{\rm stab},2}$ is needed in our proof to ensure the convergence of the discretization of the mass convection flux $\dive (\rho \bfu)$ and, to our knowledge, has never been introduced elsewhere.
Secondly, the control of the pressure in $\xLtwo(\Omega)$ relies on the stability of the discrete gradient (\ie\ the satisfaction of the so-called discrete {\em inf-sup} condition), which is not verified by colocated discretizations; note that this argument is not needed for the stability of the scheme (see the proof of {\em a priori} estimates here and \cite{gal-08-unc, eym-07-ent} for stability studies of schemes for the Navier-Stokes equations).
Assessing the effective relevance of these requirements for the discretization should deserve more work in the future.

\medskip
An easy extension of this work consists in replacing the diffusive term $-\Delta \bfu$ in \eqref{pbcont1} by its complete expression $-\mu \, \Delta \bfu -\mu/3\, \grad (\dive \bfu)$ with $\mu>0$ (\ie\ the usual form of the divergence of the shear stress tensor in a constant viscosity compressible flow).
Another less straightforward issue is the extension to more general state equations (for instance, $p=\rho^\gamma$ with $\gamma >1$); it will be the topic of a further paper.
Concerning higher order issues, let us note that the fact that the pressure is approximated by a piecewise constant function appears crucial in both stability and convergence proofs: extending this study to higher degree finite element discretizations thus certainly represents a difficult task.
Finally, let us remark that the present scheme relies on the approximation of the whole velocity vector at the interfaces. 
A less expensive scheme would be possible with a discretization $\bfu \cdot \bfn$ at the interfaces, as in the MAC scheme which is well known for the incompressible Navier-Stokes equations. 
However, such a discretization does not seem straightforward on unstructured meshes.

\bibliographystyle{plain}
\bibliography{./compiso}

\begin{thebibliography}{10}

\bibitem{boc-06-ana}
P.~Bochev, S.D. Kim, and B.-C. Shin.
\newblock Analysis and computation of least-squares methods for a compressible
  {S}tokes problem.
\newblock {\em Numerical Methods for Partial Differential Equations},
  22:867--883, 2006.

\bibitem{bre-84-sta}
F.~Brezzi and J.~Pitk{\"{a}}ranta.
\newblock On the stabilization of finite element approximations of the {S}tokes
  equations.
\newblock In W.~Hackbusch, editor, {\em Efficient Solution of Elliptic
  Systems}, volume~10 of {\em Notes Num.Fluid Mech.}, pages 11--19. Vieweg,
  1984.

\bibitem{cia-91-bas}
P.~G. Ciarlet.
\newblock Handbook of numerical analysis volume {II}~: Finite elements methods
  -- {B}asic error estimates for elliptic problems.
\newblock In P.~Ciarlet and J.L. Lions, editors, {\em Handbook of Numerical
  Analysis, Volume {II}}, pages 17--351. North Holland, 1991.

\bibitem{cro-73-con}
M.~Crouzeix and P.-A. Raviart.
\newblock Conforming and nonconforming finite element methods for solving the
  stationary {S}tokes equations {I}.
\newblock {\em Revue Fran\c{c}aise d'Automatique, Informatique et Recherche
  Op\'{e}rationnelle (R.A.I.R.O.)}, R-3:33--75, 1973.

\bibitem{ern-04-the}
A.~Ern and J.-L. Guermond.
\newblock Theory and practice of finite elements.
\newblock Number 159 in Applied Mathematical Sciences. Springer, New York,
  2004.

\bibitem{eym-00-fin}
R.~Eymard, T~Gallou{\"e}t, and R.~Herbin.
\newblock Finite volume methods.
\newblock In P.~Ciarlet and J.L. Lions, editors, {\em Handbook of Numerical
  Analysis, Volume {VII}}, pages 713--1020. North Holland, 2000.

\bibitem{eym-07-ent}
R.~Eymard and R.~Herbin.
\newblock Entropy estimate for the approximation of the compressible barotropic
  {N}avier-{S}tokes equations using a collocated finite volume scheme.
\newblock {\em in preparation}, 2007.

\bibitem{eym-06-sta}
R.~Eymard, R.~Herbin, and J.C. Latch\'e.
\newblock On a stabilized colocated finite volume scheme for the {S}tokes
  problem.
\newblock {\em Mathematical Modelling and Numerical Analysis}, 40(3):501--528,
  2006.

\bibitem{fei-04-dyn}
E.~Feireisl.
\newblock Dynamics of viscous compressible flows.
\newblock volume~26 of {\em Oxford Lecture Series in Mathematics and its
  Applications}. Oxford University Press, 2004.

\bibitem{gal-08-unc}
T.~Gallou{\"e}t, L.~Gastaldo, R.~Herbin, and J.-C. Latch{\'e}.
\newblock An unconditionnally stable pressure correction scheme for
  compressible barotropic {N}avier-{S}tokes equations.
\newblock {\em Mathematical Modelling and Numerical Analysis}, 42:303--331,
  2008.

\bibitem{gas-08-ent}
L.~Gastaldo, R.~Herbin, and J.-C. Latch{\'e}.
\newblock An entropy-preserving finite element--finite volume pressure
  correction scheme for the drift-flux model.
\newblock {\em submitted}, 2008.

\bibitem{kel-96-fin}
R.B. Kellog and B.~Liu.
\newblock A finite element method for the compressible {S}tokes equations.
\newblock {\em SIAM Journal on Numerical Analysis}, 33:780--788, 1996.

\bibitem{kel-97-pen}
R.B. Kellog and B.~Liu.
\newblock A penalized finite-element method for a compressible {S}tokes system.
\newblock {\em SIAM Journal on Numerical Analysis}, 34:1093--1105, 1997.

\bibitem{kwe-00-opt}
J.R. Kweon.
\newblock An optimal order convergence for a weak formulation of the
  compressible {S}tokes system with inflow boundary condition.
\newblock {\em Numerische Mathematik}, 86:305--318, 2000.

\bibitem{kwe-03-opt}
J.R. Kweon.
\newblock Optimal error estimate for a mixed finite element method for
  compressible {N}avier-{S}tokes system.
\newblock {\em Applied Numerical Mathematics}, 45:275--292, 2003.

\bibitem{pll-98-mat}
P.-L. Lions.
\newblock Mathematical topics in fluid mechanics -- volume 2 -- compressible
  models.
\newblock volume~10 of {\em Oxford Lecture Series in Mathematics and its
  Applications}. Oxford University Press, 1998.

\bibitem{nov-04-int}
A.~Novotn\'y and I.~Stra\v{s}kraba.
\newblock Introduction to the mathematical theory of compressible flow.
\newblock volume~27 of {\em Oxford Lecture Series in Mathematics and its
  Applications}. Oxford University Press, 2004.

\bibitem{pay-60-opt}
L.E. Payne and H.F. Weinberger.
\newblock An optimal {P}oincar\'e-inequality for convex domains.
\newblock {\em Archive for Rational Mechanics and Analysis}, 5:286--292, 1960.

\bibitem{tem-77-nav}
R.~Temam.
\newblock Navier-stokes equations.
\newblock volume~2 of {\em Studies in mathematics and its applications}. North
  Holland, 1977.

\bibitem{val-87-exi}
A.~Valli.
\newblock On the existence of stationary solutions to compressible
  {N}avier-{S}tokes equations.
\newblock {\em Annales de l'Institut Henri Poincar\'e, Section C},
  4(1):99--113, 1987.

\bibitem{ver-99-err}
R.~Verf{\"u}rth.
\newblock Error estimates for some quasi-interpolation operators.
\newblock {\em Mathematical Modelling and Numerical Analysis}, 33(4):695--713,
  1999.

\end{thebibliography}
\end{document}